\documentclass[mathpazo]{VP}
\usepackage{newpxtext} 
\usepackage{newpxmath}
\usepackage{float}
\usepackage{titlesec}
\usepackage{amsmath, amssymb, booktabs, capt-of}
\usepackage{booktabs}
\begin{document}
\title{Random Regularity  of the Vlasov-Poisson System with Random Initial Inputs in the Quasineutral Regime}
\author[W. Wang and Y. Lin]{Wenyi Wang\affil{1}\ and Yiwen Lin\affil{1}\comma\corrauth}
\address{\affilnum{1}\ School of Mathematical Sciences, Shanghai Jiao Tong University, Shanghai 200240, China}

\emails{{\tt wangwenyi2004@sjtu.edu.cn} (W. Wang), {\tt linyiwen@sjtu.edu.cn} (Y. Lin)}

\begin{abstract}
The Vlasov-Poisson system is widely used in plasma physics and other related fields. In this paper, we study the Vlasov-Poisson system with initial uncertainty in the quasineutral regime. First, we prove the uniform convergence in the Wasserstein distance between the uncertain Vlasov-Poisson system in the quasineutral regime and its quasineutral limit system with random initial inputs. This is achieved by deriving an upper bound for the Wasserstein distance and rigorously estimating each component of this bound. Furthermore, by defining a new norm with respect to the quasineutral parameter and estimating the distribution function as well as the electric field in this norm using a variable substitution, we establish the random regularity of the solutions in the quasineutral regime. This work develops a novel framework for quantifying the propagation of the initial uncertainty of the Vlasov-Poisson system in the quasineutral regime, providing a theoretical basis for designing high-performance numerical algorithms.
\end{abstract}


\maketitle


\section{Introduction}
\label{sec1}
\par In this paper, we focus on the  Vlasov-Poisson system, which has wide applications in plasma physics, semiconductor devices, optical chips, and many other fields. The initial form of the Vlasov-Poisson system is as follows:
\begin{align}
\left\{\begin{array}{c}\partial_{t} f(x,v,t)+v \cdot \nabla_{x} f+E(x,t) \cdot \nabla_{v} f=0, \\ \nabla_{x} E=\displaystyle\int_{\mathbb{R}^{d}} f d v-1, \\\left.f_{}\right|_{t=0}=f_{0}, \quad\displaystyle \int_{\mathbb{R}^{d} \times \mathbb{T}^{d}} f_{0} d v d x=1,\end{array}\right. \label{1.1}   
\end{align}
where~$\mathbb{T}=[0,1]$,~$x\in \mathbb{T}^d$~is the space variable, $t>0$ is the time variable, and $v\in \mathbb{R}^d$  is the velocity. $f=f(x,v,t)$ is the distribution function of electrons. $E=E(x,t)$ is the electric field, depending on the distribution function $f$. 
\par In 1946, Landau \cite{5} first discovered the Landau damping phenomenon, which is a common phenomenon for plasmas. Mathematically, Landau damping solution is a solution $f$ that converges to a Landau damping function $f^{\ast}$ as $t\to\infty$ in the sense of $L^{\infty}$ norm   \cite{6}. After early works conducted at the linearized level \cite{21,22}, the existence of Landau damping solution under some smoothness and smallness conditions of $f^{\ast}$ was proved in \cite{6} in 1998. Subsequently, other existence conditions were provided by Mouhot and Villani \cite{7} in 2011 and later by Bedrossian et al \cite{8}.
\par For the modeling of quasineutral regions in semiconductor devices and plasma waveguides, the Vlasov-Poisson system in the quasineutral regime is required. 
To characterize this regime, a small parameter called the quasineutral parameter $\varepsilon$ is introduced, with the quasineutral limit system derived by taking the limit $\varepsilon\to0$ \cite{4}. Due to its multiscale nature, the study in this quasineutral regime is far more challenging than that in the normal regime. In addition to the extensive development of asymptotic-preserving schemes for such multiscale problems\cite{18},
analytical results on boundedness and convergence of the solutions to the Vlasov-Poisson system in the quasineutral regime have been established in \cite{2,10,11}, which lays the foundation for this study.
\par The majority of existing studies have concentrated on the deterministic model. Nevertheless, in the engineering application of the Vlasov-Poisson system, there usually exists uncertainties due to many factors such as measurement errors and random perturbations of initial data.  Therefore a random variable $z$ is introduced to model the uncertainty. Then the initial input $f_0$ is $z$-dependent, which results in $f$ and $E$ also depending on $z$. The objective of this study is to quantify the uncertainty, that is, to investigate the propagation of the uncertainty in the initial input  and the  response uncertainty  obtained in $f$ and $E$. Conducting uncertainty quantification (UQ) analysis for the Vlasov-Poisson system is constructive in characterizing  the long-time behavior of the solutions, which is instrumental in algorithm design. In this work, we will quantify the propagation of the initial uncertainty of the Vlasov-Poisson system in the quasineutral regime through two key aspects.
\par Firstly, to quantify the distance between the Vlasov-Poisson system in the quasineutral regime and its quasineutral limit system, the Wasserstein distance is extensively used. Some properties of the Wasserstein distance are provided in \cite{12}. In 2017, Han-Kwan \cite{2} established the convergence in the Wasserstein distance between the deterministic Vlasov-Poisson system in the quasineutral regime and its quasineutral limit system based on contributions of \cite{13,14,23}. For the uncertain model, these results hold for all fixed $z$ since fixing $z$ reduces the system to the deterministic one. However, to the best of our knowledge, no prior study has addressed the uniform convergence with respect to $z$ via the Wasserstein distance between the Vlasov-Poisson system with random initial inputs in the quasineutral regime and its quasineutral limit system.
\par While the Wasserstein distance is effective for quantifying the convergence of the solutions, the random regularity analysis would provide better results by estimating both the solutions and their derivatives, consequently providing a more comprehensive characterization of the propagation of initial uncertainty. Shu and Jin \cite{3}, as well as Ding and Jin \cite{9}  achieved some results about the random regularity of $f~\text{and}~E$ in the normal regime. These random regularity results are important not only for understanding how the initial uncertainty propagates but also for performing numerical approximations. Related works can be seen in \cite{19,20}. Nevertheless, as far as we are aware, 
no existing work established the random regularity of the uncertain Vlasov-Poisson system with random initial inputs in the quasineutral regime.
\par To quantify the propagation of the initial uncertainty of the random Vlasov-Poisson system in the quasineutral regime, this work first investigates the uniform convergence of the solutions in the Wasserstein distance in the two-dimensional and three-dimensional cases. Then we define a new norm and establish the random regularity explicitly with respect to the quasineutral parameter for both the solutions and their derivatives in the one-dimensional case, which reveals the long-time behavior of the solutions. To our knowledge, this paper is a first attempt to quantify the propagation of the uncertainty of this system. The main challenges lie in the difficulty of analysis caused by uncertainty and multiscale behavior in the quasineutral regime. Firstly, from the point of the Wasserstein distance, we establish the uniform convergence in the Wasserstein distance between the uncertain  Vlasov-Poisson system with initial random inputs in the quasineutral regime and its quasineutral limit system with respect to $z$  for the two-dimensional and three-dimensional cases. This result implies that the convergence in the Wasserstein distance is insensitive to the uncertainty in the initial input, which is instructive in controlling errors in high precision numerical simulations. Building on the framework established in \cite{2} for the deterministic model, our investigation, which accounts for uncertainty in the analysis, incorporates proper constraints on $z$ and achieves uniform convergence in the Wasserstein distance by obtaining an upper bound and estimating each component of the bound. This result can be regarded as  an extension of their work to the model with uncertainty. Besides, we further investigate deeper properties of the system through the random regularity analysis. From this perspective, the random regularity of the distribution function and the electric field for the one-dimensional Vlasov-Poisson system with initial uncertainty in the quasineutral regime is studied. Inspired by the random regularity analysis in the normal regime in \cite{3},  we define a new norm with respect to the quasineutral parameter $\varepsilon$ in this work. By performing a variable substitution and estimating the solutions and their $z$-derivatives up to  total order $K$, we derive the random regularity by this norm for both the electric field $E$ and the difference between the distribution function $f$ and the Landau damping function $f^\ast$, with the convergence rate explicitly in the quasineutral parameter. 
These results contribute to the development of numerical methods such as stochastic Galerkin and stochastic collocation methods \cite{15,16,17}. The analysis on the random regularity establishes estimates not only for the solutions  but also for their derivatives, thereby revealing more intrinsic properties  of the Vlasov-Poisson system.
\par The paper is organized as follows. In Section 2 we introduce the model under different regimes. In Section 3 we review the main notations and preliminaries essential for our results.  The main results of this paper are presented in Section 4. Section 5 states some lemmas and an important proposition that support the proof of the main results. In Section 6, we prove the two main theorems of this paper. The paper is concluded in Section 7.
\section{Model}
\label{sec2}
\par In this section, we introduce the model we focus on in this work. 

\par First we present the the Vlasov-Poisson model with initial uncertainty in the normal  regime. A random variable $z$, which is in a random space $I_{z} $ with probability density function $\pi(z)$, is introduced to model the uncertainty. The initial input $f_{0}$ is $z$-dependent. Namely, $f_{0}=f_{0}(x,v,z)$. For the one-dimensional Vlasov-Poisson system, we also restrict $z$ to be one-dimensional, while in higher-dimensional cases, we do not restrict the dimension of $z$. Now, the $d$-dimensional Vlasov-Poisson system reads:
\begin{align}
\left\{\begin{array}{c}\partial_{t} f(x,v,t,z)+v \cdot \nabla_{x} f+E(x,t,z) \cdot \nabla_{v} f=0, \\ \nabla_{x} E=\displaystyle\int_{\mathbb{R}^{d}} f d v-1, \\\left.f_{}\right|_{t=0}=f_{0}, \quad\displaystyle \int_{\mathbb{R}^{d} \times \mathbb{T}^{d}} f_{0} d v d x=1.\end{array}\right.\label{2.1}
\end{align}
\par Then consider the quasineutral regime. The Debye length $\lambda _{D} $ is defined by
\begin{equation*}
 \lambda_{D}:=\left(\dfrac{\varepsilon_{0} k_{B} T_{e}}{n_{e} q_{e}^{2}}\right)^{1 / 2},   
\end{equation*}
where $\varepsilon_{0}$  is the vacuum permittivity, $k_{B}$ is the Boltzmann constant,  $T_{e}$  is the electron temperature,  $n_{e}$  is the electron density and $q_e$ is the charge of a single electron. The Debye length $\lambda _{D} $ is usually significantly smaller than the observation length $L$. Define a parameter $\varepsilon =\frac{\lambda _{D} }{L}\in(0,1) $, which is called the \textbf{quasineutral parameter}. By incorporating $\varepsilon$ into the system and applying the scaling in \cite{4}, one attains the uncertain Vlasov-Poisson system in the quasineutral regime:
\begin{align}
\left\{\begin{array}{c}\partial_{t} f_{\varepsilon}(x,v,t,z)+v \cdot \nabla_{x} f_{\varepsilon}+E_{\varepsilon}(x,t,z) \cdot \nabla_{v} f_{\varepsilon}=0, \\\varepsilon^{2} \nabla_{x} E_{\varepsilon}=\displaystyle\int_{\mathbb{R}^{d}} f_{\varepsilon} d v-1, \\\left.f_{\varepsilon}\right|_{t=0}=f_{0, \varepsilon}, \quad\displaystyle \int_{\mathbb{R}^{d} \times \mathbb{T}^{d}} f_{0, \varepsilon} d v d x=1.\end{array}\right.\label{7}
\end{align}
\par  As $\varepsilon \to 0$, the regime is called~\textbf{quasineutral limit}. The quasineutral limit of Eq. \eqref{7} is as follows:
\begin{align}
    \left\{\begin{array}{c}\partial_{t} f(x,v,t,z)+v \cdot \nabla_{x} f+E(x,t,z) \cdot \nabla_{v} f=0, \\\displaystyle\int_{\mathbb{R}^{d}} f d v=1, \\\left.f\right|_{t=0}=f_{0}, \quad \displaystyle\int_{\mathbb{R}^{d} \times \mathbb{T}^{d}} f_{0} d v d x=1.\end{array}\right.\label{8}.
\end{align}
\par Finally, as established in \cite{1}, there exists a relationship between the solutions to Eqs. \eqref{2.1} and \eqref{7}. Leveraging this connection shown in Lemma \ref{lem3.3}, we will perform a careful variable substitution that enables us to prove our main results.
\begin{lemma}\label{lem3.3}
Let $f(x,v,t,z)$ be a solution to Eq. \eqref{2.1} , then
\begin{equation}
h(x,v,t,z)=f\Bigl(\dfrac{x}{\varepsilon } ,v, \dfrac{t}{\varepsilon },z\Bigr)  \label{3.7}  
\end{equation}
is a solution to Eq. \eqref{7}.    
\end{lemma}
\par 
\begin{proof}
Let $x' = \dfrac{x}{\varepsilon}$, $t' = \dfrac{t}{\varepsilon}$. By definition, $h(x,v,t,z) = f(x',v,t',z)$. Using the chain rule,
\begin{align}
    \partial_t h(x,v,t,z)=\partial_{t'}f \dfrac{\partial t'}{\partial t}=\frac{1}{\varepsilon } \partial_{t'} f(x',v,t',z),\label{eqh1}
    \\\nabla_ x h(x,v,t,z)=\nabla_{x'}f \dfrac{\partial x'}{\partial x}=\frac{1}{\varepsilon } \nabla_ {x'}  f(x',v,t',z).\label{eqh2}
\end{align}
Since $f$ is a solution to Eq. \eqref{2.1},
\begin{gather*}
   \partial_{t'} f(x',v,t',z)+v \cdot \nabla_{x'} f+E(x',t',z) \cdot \nabla_{v} f=0, \\ \nabla_{x'} E=\displaystyle\int_{\mathbb{R}^{d}} f d v-1.
\end{gather*}
Multiplying both sides of the equation by $\dfrac{1}{\varepsilon}$, substituting \eqref{eqh1}-\eqref{eqh2} and noting that $\nabla_v h = \nabla_v f$, one has
\begin{gather}
    \partial_{t} h+v \cdot \nabla_{x} h+\frac{1}{\varepsilon}E(x',t',z) \cdot \nabla_{v} h=0.\label{2.5}
\end{gather}
For the electric field, 
since $\nabla_{x}E\Big(x',t',z\Big)=\dfrac{1}{\varepsilon}\nabla_{x'}E(x',t',z)$ by the chain rule,
one has
\begin{align}
    \nabla_x\left[\frac{1}{\varepsilon} E\left(x', t', z\right)\right]=\frac{1}{\varepsilon^2} \nabla_{x'} E=\frac{1}{\varepsilon^2}\left(\int_{\mathbb{R}^d} h d v-1\right).
\end{align}
This completes the proof that $h$ is a solution to Eq. \eqref{7}. 
\end{proof}

\par In the following part of this paper, we mainly consider Eq. \eqref{7}, which is the uncertain Vlasov-Poisson system in the quasineutral regime.

\section{Preliminaries}
\label{sec3}
\par To establish our main results, the main notations and the necessary assumptions are presented in this section.
\subsection{Notations}
{ 

\par This subsection is devoted to introducing the main notations used throughout this paper.
\par To study the uniform convergence in the Wasserstein distance, we first introduce the following notations. Let $f_\varepsilon$ denote a distribution function for the Vlasov-Poisson system Eq. \eqref{7} in the quasineutral regime. Define $\tilde{f}_\varepsilon$ by $\tilde{f}_{\varepsilon}(x, v, t,z):=f_{\varepsilon}\left(x, v-C_{\varepsilon}(x, t, z), t, z\right)$, where $C_\varepsilon$ is given in Eq. \eqref{5.2}. Denote by $g_{\varepsilon}=\int_{\mathbb{R}^{d}} \rho_{\varepsilon}^{\theta} \delta_{v=v_{\varepsilon}^{\theta}} d \mu(\theta)$ a weak solution to Eq. \eqref{7}, with $\rho_{\varepsilon}^{\theta}$ and $v_{\varepsilon}^{\theta}$ defined in Eq. \eqref{12}.  From this we further define $\tilde{g}_{\varepsilon}$ by  $\tilde{g}_{\varepsilon}:=\int_{\mathbb{R}^{d}} \rho_{\varepsilon}^{\theta} \delta_{v=v_{\varepsilon}^{\theta}+C_{\varepsilon}} d \mu(\theta)$. In addition, let $g$ represent a weak solution to the quasineutral limit system Eq. \eqref{8} and let $E_\varepsilon$ denote the electric field of the Vlasov-Poisson system Eq. \eqref{7} in the quasineutral regime. 
\par We next detail the notations employed in the random regularity analysis. Let $f$ denote the Landau damping solution to the Vlasov-Poisson system Eq. \eqref{2.1} in the normal regime. Then $h(x,v,t,z)=f(\frac{x}{\varepsilon } ,v, \frac{t}{\varepsilon },z)$ is a solution to the Vlasov-Poisson system Eq. \eqref{7} in the quasineutral regime. In this work, $f^{\ast}$ represents the Landau damping function, and $\hat{f}^{\ast}$ denotes its Fourier transform with respect to $x$ and $v$. The electric field of the Vlasov-Poisson system Eq. \eqref{2.1} in the normal regime is represented by $E$ while $E_1$ is defined as the electric field of the Vlasov-Poisson system Eq. \eqref{7} in the quasineutral regime. 
\par To describe the random regularity  of $f$ and $E$,  a norm is defined in \cite{3} as follows:
\begin{equation}\label{def_1}
    \|F\|_{\tilde{a}, t_{0}, k}=\sup  \limits_{t \geq t_{0}} t^{-k} e^{\tilde{a} t}\|F(\cdot, t)\|_{L^{\infty}},
\end{equation}
where the $L^{\infty}$ norm is taken over all variables except $t$ and $z$, $\tilde{a}>0,~t_0>0,~\text{and}~k\in \mathbb{N}$ are fixed constants.  Note that this norm Eq. \eqref{def_1} is defined for the random regularity analysis in the normal regime, while our study focuses on the quasineutral regime, which requires the construction of a new norm with respect to the quasineutral parameter $\varepsilon$.
\par To attain the random regularity in the quasineutral regime, which is full of challenges due to its multiscale nature, we define a new norm with respect to the quasineutral parameter $\varepsilon$:
\begin{equation}
    \left \| F \right \|  _{a,t_{0},k,m } =\displaystyle\sup \limits_{t\ge t_{0} } t^{-k} \exp\biggl[{\Bigl(a+\frac{1}{\varepsilon ^{m} }\Bigr)t }\biggr] \left \| F(\cdot ,t) \right \| _{L^{\infty }}  ,\label{3.2}
\end{equation}
where the $L^{\infty}$ norm is taken over all variables except $t$ and $z$, $a\in \mathbb{R},~t_0>0,~k\in \mathbb{N},~\text{and}~m\in \mathbb{N_+}$ are fixed constants. Such norm is critical for our result.
\endgraf  
}
\subsection{Basic assumptions in the Wasserstein distance}
\par Let us introduce the Wasserstein distance to quantify the discrepancy between the Vlasov-Poisson system in the quasineutral regime and its quasineutral limit system.
\par \begin{definition}[Wasserstein distance]
  Let  $ \mathcal{P}_{q}(\mathbb{T}^{d}\times\mathbb{R}^{d})$  be the collection of all probability measures $ \gamma $ on  $\mathbb{T}^{d}\times\mathbb{R}^{d} $ with finite  $q$  moment: for some $ x_{0} \in {\mathbb{T}^{d}\times\mathbb{R}^{d}} , \int_{\mathbb{T}^{d}\times\mathbb{R}^{d}} d(x, x_{0})^{q} \mathrm{d} \gamma(x)<+\infty$~, where~$d(\cdot , \cdot )$ is the Euclidean distance.
Then the  $q$-Wasserstein distance between two probability measures  $ \nu_1 $  and  $ \nu_2 $  in $ \mathcal{P}_{q}(\mathbb{T}^{d}\times\mathbb{R}^{d}) $ is defined as:
\begin{equation*}   
W_{q}(\nu_1, \nu_2):=\displaystyle\left(\inf \int_{(\mathbb{T}^{d} \times \mathbb{R}^{d})^2} d(x, y)^{q} \mathrm{~d} \gamma(x, y)\right)^{1 / q},
\end{equation*}
where the infimum is taken over all measures $\gamma$  on  $(\mathbb{T}^{d} \times \mathbb{R}^{d})^2 $ with marginals $ \nu_1 $ and $ \nu_2 $ on the first and second factors respectively. 
\end{definition}
\par In order to obtain a solution to the Vlasov-Poisson system Eq. \eqref{7} in the quasineutral regime, which is essential for proving our main results in the Wasserstein distance, the multi-fluid pressureless Euler-Poisson system is now introduced: 
\begin{align}
    \left\{\begin{array}{c}\partial_{t} \rho_{\varepsilon}^{\theta}(x,t,z)+\nabla_{x} \cdot\left(\rho_{\varepsilon}^{\theta} v_{\varepsilon}^{\theta}(x,t,z)\right)=0 , \\\partial_{t} v_{\varepsilon}^{\theta}+v_{\varepsilon}^{\theta} \cdot \nabla_{x} v_{\varepsilon}^{\theta}=E_{\varepsilon} , \\\varepsilon^{2} \nabla_{x} E_{\varepsilon}=\displaystyle\int_{\mathbb{R}^d} \rho_{\varepsilon}^{\theta} d \mu(\theta)-1 , \\\left.\rho_{\varepsilon}^{\theta}\right|_{t=0}=\rho_{0, \varepsilon}^{\theta}, \left.v_{\varepsilon}^{\theta}\right|_{t=0}=v_{0, \varepsilon}^{\theta}, \end{array}\right.\label{12}
\end{align}
 where 
 $d \mu(\theta)=c_{d} \frac{d \theta}{1+|\theta|^{d+1}}$ with $c_d$ a constant only dependent on $d$.
Specifically, let $\rho^{\theta}_\varepsilon $ be a solution to Eq. \eqref{12}, then 
 \begin{equation*}
  \displaystyle g_{\varepsilon}(x,v,t,z)=\int_{\mathbb{R}^{d} } \rho_{\varepsilon}^{\theta}(x,t,z)\delta_{v=v_{\varepsilon}^{\theta}(x,t,z)} d \mu(\theta)   
 \end{equation*}
 is a weak solution to Eq. \eqref{7}.
\par The formal limit system of Eq. \eqref{12} is the multi-fluid 
 incompressible Euler system:
\begin{align}
    \left\{\begin{array}{c}\partial_{t} \rho^{\theta}(x,t,z)+\nabla_{x} \cdot\left(\rho^{\theta} v^{\theta}(x,t,z)\right)=0, \\\partial_{t} v^{\theta}+v^{\theta} \cdot \nabla_{x} v^{\theta}=E , \\\operatorname{curl} E=0, ~\displaystyle\int_{\mathbb{T}^{d}} E d x=0, \\\displaystyle \int_{\mathbb{R}^d} \rho^{\theta} d \mu(\theta)=1 , \\\left.\rho^{\theta}\right|_{t=0}=\rho_{0}^{\theta}, \left.v^{\theta}\right|_{t=0}=v_{0}^{\theta}, \end{array}\right.\label{4.2}
\end{align}
 where 
 $ \lim\limits_{\varepsilon  \to 0} \rho _{0, \varepsilon }^{\theta } =\rho _{0}^{\theta } \text { (suppose the limit exists), }  ~\text{and}~v_{0, \varepsilon}^{\theta}=v_{0}^{\theta}=\theta \text {. }$
 \par Let $\rho^{\theta}$ denote a solution to Eq. \eqref{4.2}, then 
 \begin{equation*}
     \displaystyle g(x,v,t,z)=\int_{\mathbb{R}^{d} } \rho^{\theta}(x,t,z)\delta_{v=v^{\theta}(x,t,z)} d \mu(\theta) 
 \end{equation*}
 defines a weak solution to Eq. \eqref{8}.
 \par  To establish the uniform convergence in the Wasserstein distance, we introduce the following assumptions based on \cite{2}.
 Specifically, (i)-(vi) concern the initial distributions and initial electric field, while (vii) and (viii) relate to the initial Wasserstein distance, as summarized in the table below. 
 
 \begin{assumption}[Assumptions for initial distributions and initial Wasserstein distance]\label{ass3.2}
 \begin{center}
  \label{tab:assump3.1}
  \renewcommand{\arraystretch}{1.3}
  \begin{tabular}{lp{0.9\linewidth}}
  & Suppose that the following assumptions hold for all $z\in I_{z}$.\\
    \midrule
    & The initial distributions and initial electric field satisfies the following (i)-(vi):\\
    (i) & $f_{\varepsilon }| _{t=0} : = f_{0,\varepsilon} = 0$ if $|v| > \frac{1}{\varepsilon^\gamma}$, where $\gamma > 0$ is a fixed constant (initial truncation in velocity space). \\
    (ii) & 
      $f_{0,\varepsilon} = g_{0,\varepsilon} + u_{0,\varepsilon}$, where $g_{0,\varepsilon}$ is a sequence of continuous functions with uniform compact support in $v$, satisfying: \\
      & \quad (1) $\sup\limits_{\varepsilon \in (0,1)} \sup\limits_{v \in \mathbb{R}^d} \left(1 + |v|^2\right) \| g_{0,\varepsilon}(\cdot, v,z) \|_{B_{\delta_0}} \leq C_0$; \\
      & \quad (2) $\left\| \int_{\mathbb{R}^d} g_{0,\varepsilon}(\cdot, v,z) \, dv - 1 \right\|_{B_{\delta_0}} \leq C_0 \varepsilon$,  $\forall \varepsilon \in (0,1)$,\\
      & where $\| g \|_{B_\delta} := \sum_{k \in \mathbb{Z}} | \widehat{g}(k) | \delta^{|k|}$ with $\widehat{g}(k)$ the $k$-th Fourier coefficient of $g$ in $x$, and $\delta_0 > 0, C_0 > 1$ are fixed numbers. \\
    (iii) & As $\varepsilon \to 0^+$, $g_{0,\varepsilon} \to g_0$ in the sense of distribution, and $g(x, v, 0,z) = g_0$.\\
    (iv)& $u_{0, \varepsilon}$ is a sequence of functions satisfying for all  $\varepsilon>0$, $W_{2}\left(f_{0, \varepsilon}, g_{0, \varepsilon}\right)=\psi(\varepsilon,z) .$\\ 
   (v) & $ \| f_{0, \varepsilon }   \| _\infty\le C_{0} , ~ \frac{1}{2} \int_{\mathbb{R}^{d} \times \mathbb{T}^{d}} f_{0, \varepsilon }|v|^{2} dvd x\le C_{0}$, where $ \| \cdot  \| _\infty$ is the $L_\infty$ norm over $x$. \\
   (vi)& $\parallel E_{\varepsilon } (x,0,z) \parallel _{\infty } \le C_{0}$.\\
   \midrule
    &The initial Wasserstein distance $\psi(\varepsilon,z) = W_2(f_{0,\varepsilon}, g_{0,\varepsilon})$ satisfies the following (vii) and (viii): \\
    (vii)&$\psi(\varepsilon,z) \leq d$, where $d = 2, 3$ is the system dimension. \\
    (viii) & For a fixed $T > 0$, $R_\varepsilon\left(T, \psi(\varepsilon,z),z\right) \leq d$, where \\
      & \quad $\displaystyle R_\varepsilon(t, x,z) := 16d \exp\left\{ \log\left(\frac{x}{16d}\right) \exp\left[ C_0 \int_0^t \tilde{A}(s,z) \, ds \right] \right\}$, \\
      & \quad $ \tilde{A} (t,z):=  1+\dfrac{\sqrt{\left\|\rho_{g_\varepsilon}(t,z)\right\|_{\infty}}\biggl[\max \left\{\left\|\rho_{f_\varepsilon}(t,z)\right\|_{\infty}, \left\|\rho_{g_\varepsilon}(t,z)\right\|_{\infty}\right\}\biggr]^{1 / 2}+\left\|\rho_{f_\varepsilon}(t,z)-1\right\|_{\infty} }{\varepsilon^{2}} $, \\
      & with $\rho_{f_\varepsilon} = \int_{\mathbb{R}^d} f_\varepsilon \, dv$, $\rho_{g_\varepsilon} = \int_{\mathbb{R}^d} g_\varepsilon \, dv$, and $g_\varepsilon(x, v, t,z) = \int_{\mathbb{R}^d} \rho_{\varepsilon}^\theta(x, t,z) \delta_{v= v_{\varepsilon}^\theta(x, t,z)} \, d\mu(\theta)$ a weak solution to Eq. (2.2) with $g_{\varepsilon }| _{t=0} = g_{0,\varepsilon}$.\\
    \bottomrule
  \end{tabular}
  \end{center}
\end{assumption}
\subsection{Basic assumptions on the regularity}
Now we introduce two assumptions for the random regularity analysis: Assumption \ref{ass3.4} concerning the Landau damping function from \cite{6}, and Assumption \ref{ass3.5} regarding the parameter constraints from \cite{3}.
\par If a solution to Eq. \eqref{2.1} satisfies
\begin{equation}
\left \| f(x,v,t,z)- f^{\ast } (x-vt,v,z) \right \| _{L^{\infty }(x,v)}\to 0 ~\text{as}~t\to \infty ,  \label{17*}
\end{equation}
where $f^{\ast }$ is a known function called \textbf{Landau damping}, then such a solution is called \textbf{Landau damping solution}. If Assumption \ref{ass3.4} below holds, then there exists an initial input at $t=t_0$ such that Eq. \eqref{2.1} admits a Landau damping function.
\begin{assumption}[Assumptions on Landau damping function]\label{ass3.4}
Suppose the Landau damping function satisfies:
\par$  (H1)~ \hat{f}^{*}(k_{x}, k_{v}) \leq \frac{a_{1}}{1+k_{x}^{2}} e^{-\tilde{a}|k_{v}|} ~\text{(Smoothness)},$ where 
\begin{equation*}
\displaystyle\hat{f}^{*}\left(k_{x}, k_{v}\right)=\frac{1}{2\pi} \iint_{\mathbb{R}\times\mathbb{T}} f^{*}(x, v) e^{i\left(k_{x} x+k_{v} v\right)} dvdx.    
\end{equation*}
\par $ (H2)~|f^{*}(x, v)| \leq \frac{a_{2}}{1+v^{4}} ~\text{(Decay).}$
\par ($H3$)~There exist constants such that~$\tilde{a} \geq 15 \sqrt{a_{2}}, $  and $t_{0} \geq \max \{0, \frac{1}{\tilde{a}} \log (8 a_{1})\}.$    
\end{assumption}

\begin{assumption}[Assumptions on parameters]\label{ass3.5} Assume that the parameters satisfy the following inequalities:
 \par $(A1)~  \tilde{a} \geq \max \{1,15 \sqrt{a_{2}}\} .$
\par$(A2)~  t_{0} \geq \max \{2,4 K, \frac{1}{\tilde{a}} \log (8 a_{1})\}.$
\par$(A3) ~ \frac{50 C_{E}}{\tilde{a}}(3 / \tilde{a})^{3} e^{-3} \leq 1 ,~\text{i.e.} \: \frac{50 C_{E}}{\tilde{a}} t_{0}^{3} e^{-\tilde{a} t_{0}} \leq 1,~\text{where} ~C_{E}=\frac{240 a_{1} a_{2}}{\tilde{a}}+4 a_{1}.$
\par$(A4)~  8 e \leq \frac{1}{20 a_{2}} .$
\par$(A5) ~ 8 C_{E} \leq \tilde{a}^{2} .$   
\end{assumption}

\section{Main results}
\label{sec4}
\par 
In this paper, we mainly achieve results in the following two perspectives:
\par(A) We prove the uniform convergence with respect to $z$ in the Wasserstein distance between the Vlasov-Poisson system Eq. \eqref{7} with initial uncertainty in the quasineutral regime and its quasineutral limit system Eq. \eqref{8} in the two-dimensional and three-dimensional cases. This result is established by imposing fitting restrictions on $z$, attaining an upper bound for this distance and then estimating each part of this upper bound.
\par(B) By defining a new norm with respect to $\varepsilon$ and performing a variable substitution to estimate the solutions to Eq. \eqref{7} in the one-dimensional case as well as their derivatives under this norm, we establish the random regularity of the Vlasov-Poisson system with initial uncertainty in the quasineutral regime.
\par For higher-dimensional cases ($d=2, 3$), we establish the convergence of the distribution function of the uncertain Vlasov-Poisson system in the quasineutral regime by imposing the boundedness of the density, velocity, and their respective errors  with respect to the limit counterparts in Theorem \ref{thm1}. For $d=1$, by defining a norm with respect to the quasineutral parameter, we establish the random regularity explicitly with respect to the quasineutral parameter $\varepsilon$ for both the solutions and their derivatives in Theorem \ref{thm2}, which indicates insensitivity to the propagation of uncertainty in the initial data and suggests that the random regularity can persist in long-time in the quasineutral regime, revealing deeper properties of the system. Note that for $d=1$, analogous convergence analysis in the Wasserstein distance can also be conducted. However, in this work, we establish a more profound  random regularity analysis for the one-dimensional scenario. On the other hand, due to additional challenges in higher-dimensional cases, our random regularity analysis is currently restricted to one dimension. The extension of random regularity analysis to higher dimensions are reserved for future investigations.
\par To address the additional challenges posed by the randomness and the quasineutral regime, we construct a set $A_{z,\delta,M,z_{0}}$ and a new norm $\left \| \cdot \right \|  _{a,t_{0},k,m }$ respectively. The set $A_{z,\delta,M,z_{0}}$ plays a vital role in achieving uniform convergence in the Wasserstein distance for the random model, and the new norm is essential for proving the random regularity in the quasineutral regime.

\subsection{The convergence in the Wasserstein distance for the uncertain model}
\par 
To address the challenges posed by randomness, we construct a critical set $A_{z, \delta,M,z_0}$ and establish the uniform convergence in the Wasserstein distance with respect to $z$ on it.

For $d=2,3$, Let $f_\varepsilon(\cdot,z)$ and $E_\varepsilon(\cdot,z)$ be a pair of weak solutions to Eq. \eqref{7}. We define 
\begin{gather}
  A_{z, \delta,M,z_0} =\left \{ z|~\forall 0<\varepsilon <\delta ,\left \| \rho _{\varepsilon }^{\theta }(\cdot,z) \right \|_{\infty} \le M\left \| \rho _{\varepsilon }^{\theta } (\cdot,z_0)\right \|_{\infty}, \right.\nonumber\\\left.\left \| v^{\theta }(\cdot,z)\right \|_{\infty} \le M\left \| v^{\theta }(\cdot,z_0)\right \|_{\infty} , \right.\nonumber\\\left.\left \|\rho _{\varepsilon}^{\theta }(\cdot,z) -\rho ^{\theta }(\cdot,z) \right \|_{\infty} \le M\left \|\rho _{\varepsilon }^{\theta }(\cdot,z_0) -\rho ^{\theta }(\cdot,z_0) \right \|_{\infty} , \right.\nonumber\\\left.\left \| v _{\varepsilon}^{\theta } (\cdot,z)-v ^{\theta }(\cdot,z)+C_{\varepsilon}(\cdot,z)   \right \|_{\infty}\le M\left \| v _{\varepsilon }^{\theta } (\cdot,z_0)-v ^{\theta }(\cdot,z_0)+C_{\varepsilon}(\cdot,z_0)   \right \|_{\infty} \right \}\text{,}
\end{gather}
where $0<\delta<1~\text{and}~M>1$ are fixed $z$-independent constants; $z_0$ is a fixed $z$-independent vector in $I_z$, and
\begin{equation}
 C_{\varepsilon}(x,t,z):=-\dfrac{1}{i}\left(d_{+}(x,t,z) \exp\Bigl({\dfrac{i t}{\varepsilon}}\Bigl)-d_{-}(x,t,z) \exp\Bigl({-\dfrac{i t}{\varepsilon}}\Bigl)\right), \label{5.2}
\end{equation}
with $d_{\pm } $  the solutions to:
 \begin{align*}
 \operatorname{curl} d_{ \pm}&=0, \quad \operatorname{div}\Biggl(\partial_{t} d_{ \pm}+\left(\displaystyle\int_{\mathbb{R}^d}\rho^{\theta} v^{\theta} d\mu(\theta) \cdot \nabla\right) d_{ \pm}\Biggl)=0, \\\text{and}\  \operatorname{div} d_{ \pm}(x, 0,z)&=\lim\limits _{\varepsilon \rightarrow 0} \operatorname{div} \dfrac{\varepsilon E_{\varepsilon}(x, 0,z) \pm i j^{\varepsilon}(x, 0,z)}{2}~\bigg(j^{\varepsilon}:=\int_{\mathbb{R}^d} \rho_{\varepsilon}^{\theta} v_{\varepsilon}^{\theta} d\mu  (\theta)\bigg)\text {. }
 \end{align*}
In practice, one may first fix a sample $z_0$, then there exist appropriate $\delta>0$ and $M>1$.  Subsequently, one can find $z$ that satisfies the conditions to define the set $A_{z,\delta,M,z_0}$. The construction of $A_{z,\delta,M,z_0}$ is essential to establishing the uniform convergence results.
\par Theorem \ref{thm1} below gives the uniform convergence with respect to $z$ in the Wasserstein distance between the uncertain Vlasov-Poisson system in the quasineutral regime and its quasineutral limit system  in the two-dimensional ($d=2$) and three-dimensional ($d=3$) cases:
\begin{theorem}\label{thm1}
For $d=2,3$, Let $f_\varepsilon(\cdot,z)$ and $E_\varepsilon(\cdot,z)$ be a pair of weak solutions to Eq. \eqref{7}, $g(\cdot,z)$ be a weak solution to Eq. \eqref{8}. Suppose that Assumption \ref{ass3.2} (i)-(vi) holds and
\begin{equation}
 \displaystyle\int_{\mathbb{R}^d\times\mathbb{T}^d}\Big(|v|^{2}+U_{\varepsilon}(x, 0,z)\Big) f_{\varepsilon}( x, v, 0,z) d v d x\leq C_{0}, \label{22}   
\end{equation}
where $ U_\varepsilon(x,t,z)=\int_{0}^{x} E_\varepsilon(s,t,z)ds$. Then there exists an initial Wasserstein distance $\psi (\varepsilon,z )=W_{2} (f_{0, \varepsilon },g_{0, \varepsilon })$ such that $\sup \limits _{t \in[0, T]} W_{1}(\tilde{f_{\varepsilon}}(t,z), g(t,z))$ uniformly converges to 0 in~$A_{z, \delta,M,z_0}$ as $\varepsilon \to 0$, where
\begin{equation}
 \displaystyle\tilde{f_{\varepsilon}}(x, v,t,z):=f_{\varepsilon}\Big(x, v-C_{\varepsilon}(x,t,z),t,z\Big).   
\end{equation}

\begin{remark}
          The intuitive explanation for constructing the set $A_{z,\delta,M,z_0}$ is to ensure that the norms of the density $\rho _{\varepsilon }^{\theta }(\cdot,z)$ , the velocity $v^{\theta }(\cdot,z)$, the density convergence error $\rho _{\varepsilon}^{\theta }(\cdot,z) -\rho ^{\theta }(\cdot,z)$, and the velocity convergence error $v _{\varepsilon}^{\theta } (\cdot,z)-v ^{\theta }(\cdot,z)+C_{\varepsilon}(\cdot,z)$ are each bounded by $M$ times (with $M>1$) corresponding norms at point $z_0$, i.e.,  those of $\rho _{\varepsilon }^{\theta }(\cdot,z_0),v^{\theta }(\cdot,z_0),\rho _{\varepsilon}^{\theta }(\cdot,z_0) -\rho ^{\theta }(\cdot,z_0), \text{and} ~v _{\varepsilon}^{\theta } (\cdot,z_0)-v ^{\theta }(\cdot,z_0)+C_{\varepsilon}(\cdot,z_0)$, respectively. Such $L^\infty$ bounds on the solution ensure the boundedness of $\sup\limits_{t\in[0,T]} W_{1}\left(\tilde{g}_{\varepsilon}(t,z), g\right(t,z))$ (see Section 6.1 for a detailed proof), which is critical for the first step of the proof of Theorem 4.1.  Note that the set $A_{z,\delta,M,z_0}$, defined for any arbitrary point $z_0\in I_z$, is always non-empty set since $M > 1$ ensures $z_0\in A_{z,\delta,M,z_0}$. 
\end{remark}

\end{theorem}
\par Theorem \ref{thm1} reveals that the convergence in the Wasserstein distance is insensitive to the uncertainty in the initial data. 
\subsection{The random regularity in the quasineutral regime}
\par In this subsection, we investigate the random regularity by estimating derivatives in space and the random variable, thereby uncovering more profound properties of solutions in the one-dimensional case. 
\par Recalling that the new norm \eqref{3.2} with respect to the quasineutral parameter $\varepsilon$ is defined as
\begin{equation*}
    \left \| F \right \|  _{a,t_{0},k,m } =\displaystyle\sup \limits_{t\ge t_{0} } t^{-k} \exp\biggl[{\Bigl(a+\frac{1}{\varepsilon ^{m} }\Bigr)t }\biggr] \left \| F(\cdot ,t) \right \| _{L^{\infty }}  ,
\end{equation*}
with 

$a\in \mathbb{R},~t_0>0,~k\in \mathbb{N},~\text{and}~m\in \mathbb{N_+}$.
\par This norm plays a crucial role in establishing the random regularity of the  Vlasov-Poisson system in the quasineutral regime. By this norm, our main resultabout the random regularity of the one-dimensional Vlasov-Poisson system in the quasineutral regime is stated in Theorem \ref{thm2}:

\begin{theorem}\label{thm2}
    For $d=1$, Suppose that there is a fixed $\delta_1\in(0,1)$ such that for each fixed $\varepsilon\in(0, \delta_1)$, there exist~$\tilde{a} (\varepsilon )=a+\frac{1}{\varepsilon^{m}},~a_{1}(\varepsilon ),~a_{2}(\varepsilon ),~\varepsilon\text{-independent } t_{0}, \varepsilon$\mbox{-independent} $f^{\ast }$ and all of its $x,v,z$-derivatives up to total order $K$ that satisfy Assumptions \ref{ass3.4}-\ref{ass3.5}, where $K$ is a fixed non-negative integer. Let $f(x,v,t,z)$ be the Landau damping solution to Eq. \eqref{2.1} that satisfies Eq. \eqref{17*}. Then $h(x,v,t,z)=f(\frac{x}{\varepsilon } ,v, \frac{t}{\varepsilon },z)$ is a solution to Eq. \eqref{7}. Let $E_1$ be the electric field corresponding to $h$, i.e., $h~\text{and}~E_1$ form a pair of solutions to Eq. \eqref{7}.  Then the following convergence holds:
\begin{equation}
    \left\|\partial_{z}^{k}\left[h(x, v, t, z)-f^{*}\Big(\dfrac{x-v t }{\varepsilon } , v, z\Big)\right]\right\|_{a, t_{0}, 1,m} \to 0,  ~\text{as}~\varepsilon \to 0, ~\text{for }0\le k \le K.\label{17}
\end{equation}
The convergence rate reaches 
\begin{equation}
 \displaystyle O\Biggl(\dfrac{1}{\varepsilon } \exp{\biggl[- \Bigl(\frac{1}{\varepsilon ^{m+1} } -\frac{1}{\varepsilon ^{m}}+\frac{a}{\varepsilon } \Bigr)t_{0}\biggl ] }\Biggr). \label{4.7}  
\end{equation}
$\text { Furthermore, for all }~l \geq 1,~k\geq0,~l+k \leq K,~\text{except}~ k=0,~l=1 ,$ we have
\begin{equation}
    \left \| \partial _{x}^{l} \partial _ {z}^{k} E_{1}  \right \| _{a,t_{0},l,m} \to 0 ,  ~\text{as}~\varepsilon \to 0,\label{18}
\end{equation}
and when $l=0$, for all $1\le k\le K$ and any $x_0\in[0,1]$, 
\begin{equation}
 \left \| \partial _ {z}^{k} \big[E_{1}(x,t,z)-E_{1}(x_0,t,z)\big]  \right \| _{a,t_{0},0,m} \to 0 ,  ~\text{as}~\varepsilon \to 0.\label{18*}   
\end{equation}
Now the convergence rate is
\begin{equation}
 \displaystyle O\Biggl(\dfrac{1}{\varepsilon^{2l+1} } \exp{\biggl[- \Bigl(\frac{1}{\varepsilon ^{m+1} } -\frac{1}{\varepsilon ^{m}}+\frac{a}{\varepsilon } \Bigr)t_{0}\biggl ] }\Biggr).  \label{4.10}
\end{equation}
\begin{remark}
The convergence orders Eqs. \eqref{4.7} and \eqref{4.10}  both exceed $O(\exp(-\varepsilon^{-m}))$ as $\varepsilon \to 0$. As $m$ increases, the value of $\parallel \cdot\parallel _{a,t_{0},k,m }$ for the same function increases while higher degree of convergence is obtained.  For $\varepsilon=1$, these results are consistent with the random regularity results for the uncertain Vlasov-Poisson system in the normal regime established in \cite{18}. 
\end{remark}
\begin{remark}
 Extending Theorem \ref{thm2} to higher dimensions is non-trivial as higher-dimensi\-onal scenarios introduce additional analytical challenges. For instance, $\nabla^{l}_{x}$ becomes an $l$-th order tensor in higher-dimensional cases, and integrating such derivative terms proves more intricate than in the one-dimensional case, thereby making it difficult to obtain analogous estimates to Eqs. \eqref{6.38} and \eqref{6.18} in the one-dimensional case. Therefore, the investigation of random regularity in higher dimensions is reserved for future work. 
\end{remark}
\end{theorem} 

\par Theorem \ref{thm2} characterizes the long-time behavior of the solutions. Specifically, Eq. \eqref{17} demonstrates that  when time is large, the $z$-dependence of the particle distribution is dominated by the $z$-dependence of the time-asymptotic profile and insensitive to the random perturbation in the initial input. Eqs. \eqref{18} and \eqref{18*} indicate that when time is large, the electric field is insensitive to the random perturbation in the initial input. More precisely, as time approaches infinity, the particle distribution converges to the time-asymptotic profile at an exponential rate, while the fluctuation of the electric field decays exponentially in time. Moreover, as $\varepsilon\to0$, the convergence rate increases exponentially.
\section{Necessary proposition and lemmas}
\label{sec5}
\par In this section, we cite some lemmas for 
 the proof of Theorem \ref{thm1} and an important proposition that will be used in the proof of Theorem \ref{thm2}.
\par The proof of Theorem \ref{thm1} is based on the following lemmas  from \cite{2}, namely Lemmas \ref{lem5.1}-\ref{lem5.4}, which estimate the Wasserstein distance for the deterministic model in the two-dimensional ($d=2$) and three-dimensional ($d=3$) cases.
\par First, we refer to an upper bound of $\sup \limits _{t \in[0, T]} W_{1}(\tilde{f_{\varepsilon}}(t), g(t))$ for the deterministic model provided in \cite{2}  as it contributes to estimating $\sup \limits _{t \in[0, T]} W_{1}(\tilde{f_{\varepsilon}}(t,z), g(t,z))$ for the uncertain model:
\begin{lemma}[Estimate of the $W_{1}$ distance]\label{lem5.1}
 Under Assumption \ref{ass3.2} (i)-(vi), 
 the following estimates hold:
\begin{equation*}
 W_{1}\left(\tilde{f_{\varepsilon}}, g\right) \leq W_{1}\left(\tilde{f_{\varepsilon}}, \tilde{g}_{\varepsilon}\right)+W_{1}\big(\tilde{g}_{\varepsilon}, g\big)~\text{(triangle inequality),}   
\end{equation*} 
and
\begin{equation*}
 W_{1}\left(\tilde{f}_{\varepsilon}, \tilde{g}_{\varepsilon}\right) \le\left(1+C_T\right) W_{1}\left(f_{\varepsilon}, g_{\varepsilon}\right)\le \left(1+C_T\right) W_{2}\left(f_{\varepsilon}, g_{\varepsilon}\right),   
\end{equation*}
 where 
 \begin{align*}
   \displaystyle\tilde{g}_{\varepsilon}(x, v,t)=&\int_{\mathbb{R}^{d}} \rho_{\varepsilon}^{\theta}(x, t) \delta_{v=v_{\varepsilon}^{\theta}(x, t)+C_{\varepsilon}(x, t)} d \mu(\theta),
\end{align*}
 with $C_{T}$ a constant only dependent on $T$.
\end{lemma}
\begin{proof}
The first inequality in Lemma \ref{lem5.1} follows from the triangle inequality. By the Kantorovich duality, one has 
\begin{equation*}
    W_{1}(\tilde{f}_\varepsilon, \tilde{g}_\varepsilon)=\sup _{\|\varphi\|_{L i p} \leq 1}\left[\int \Big(\varphi(x, v+C_\varepsilon) f_\varepsilon(x,v)-\varphi(x, v+C_\varepsilon)g_\varepsilon(x,v)\Big)dvdx\right].
\end{equation*}
Denote $\tilde{\psi} (x, v)=\varphi(x, v+C_\varepsilon(\cdot,x)$, computing the gradient of $\tilde{\psi}$ we derive
\begin{equation}
    \|\tilde{\varphi} \|_{L i p} \leq\left(1+\left\|D_{x} C_\varepsilon\right\|_{L^{\infty}}\right)\|\varphi\|_{L i p} \leq1+\left\|D_{x} C_\varepsilon\right\|_{L^{\infty}}\le1+C_T.\label{5.1}
\end{equation}
The last inequality holds from the definition of $C_\varepsilon$ and Assumption \ref{ass3.2} (i)-(vi).
From Eq. \eqref{5.1} we obtain the second inequality in Lemma \ref{lem5.1}.
\end{proof}
\par For the deterministic model \cite{2}, an upper bound of $\sup\limits_{t\in[0,T]} W_{1}(\tilde{g}_{\varepsilon}, g)$ and its convergence as $\varepsilon\to0$ are stated in Lemma \ref{lem5.2}, with Lemma \ref{lem5.3} giving an estimate of $\sup \limits _{t \in[0, T]} W_{2}(f_{\varepsilon}, g_{\varepsilon}) $, both of which are of assistance to our proof. 
 \begin{lemma}\label{lem5.2}
 Under  Assumption \ref{ass3.2} (i)-(vi), as $\varepsilon\to0$, 
\begin{align*}
   \sup\limits_{t\in[0,T]} W_{1}\left(\tilde{g}_{\varepsilon}, g\right) & \leq \sup\limits_{t\in[0,T]}\Biggl\{ \sup \limits _{\theta \in \mathbb{R}^d}\left\|\rho_{\varepsilon}^{\theta}(x,t)\right\|_{\infty}\int_{\mathbb{R}^d}\left\|v_{\varepsilon}^{\theta}(x, t)+C_{\varepsilon}(x, t)-v^{\theta}(x, t)\right\|_{\infty} d \mu(\theta) \nonumber\\&+ \int_{\mathbb{R}^d}\left\|\rho_{\varepsilon}^{\theta}(x,t)-\rho^{\theta}(x,t)\right\|_{\infty} d \mu(\theta)\left(1 / 2+\sup \limits _{\theta \in \mathbb{R}^d}\left\|v^{\theta}(x, t)\right\|_{\infty}\right)\Biggr\}\to 0.   
   \end{align*}
\end{lemma}

\begin{lemma}[Estimate of the $W_2$ distance]\label{lem5.3}
If $\psi (\varepsilon )=W_{2} (f_{0, \varepsilon },g_{0, \varepsilon })$  satisfies (vii) and (viii) in Assumption \ref{ass3.2}, one has the following upper bound:
\begin{equation}
    \sup \limits _{t \in[0, T]} W_{2}\left(f_{\varepsilon}, g_{\varepsilon}\right) \leq 16 d \exp \Biggl\{\operatorname { l o g } \biggl( \dfrac { \psi ( \varepsilon ) } { 1 6 d } \biggr ) \operatorname {exp} \left[C _ { 0 } T \dfrac { 1 } { \varepsilon ^ { 2 } } \left(1+\left\|\rho_{f_{\varepsilon}}\right\|_{L^{\infty}\left([0, T] ; L_{x}^{\infty}\right)}+\left\|\rho_{g_{\varepsilon}}\right\|_{L^{\infty}\left([0, T] ; L_{x}^{\infty}\right)}\right)\right]\Biggr\}.\label{41}
\end{equation}    
\end{lemma} 
 \par The following estimates in Lemma \ref{lem5.4} will be helpful in 
 estimating the right-hand side of Eq. \eqref{41}, where we give the 
 estimates of $\|\rho_{f_\varepsilon}\|_{\infty}~\text{and}~\|\rho_{g_\varepsilon}\|_{\infty}$ as in \cite{2} for the convenience of the subsequent proof of Theorem \ref{thm1}.
\begin{lemma}[Estimates of $\|\rho_{f_\varepsilon}\|_{\infty}~\text{and}~\|\rho_{g_\varepsilon}\|_{\infty}$]\label{lem5.4}
 Suppose that Assumption \ref{ass3.2} (i)-(vi) and
\begin{equation*}
 \int_{\mathbb{R}^d\times\mathbb{T}^d}\Big(|v|^{2}+U_{\varepsilon}(x, 0)\Big) f_{\varepsilon}( x, v, 0) d v d x\leq C_{0}   
\end{equation*}
hold. Let $T  > 0$ be fixed. Then 
\par(1)(i) In the two-dimensional case, for~$\forall\beta   > 2$, there exists a constant $C _{\beta} > 0$ s.t. for all $\varepsilon \in (0,1)$ and all $t\in [0,T]$,
\begin{equation*}
    \displaystyle\left\|\rho_{f_\varepsilon}\right\|_{\infty} \leq C_{\beta} / \varepsilon^{2 \max \{\beta, \gamma\}}.
\end{equation*}
\par~~~(ii) In the three-dimensional case, there is a constant  $C_1 > 0$ s.t. for all $\varepsilon \in (0,1)$ and all $t\in [0,T]$,
\begin{equation*}
    \left\|\rho_{f_\varepsilon}\right\|_{\infty} \leq C_1 / \varepsilon^{ \max \{38, 3\gamma\}}.
\end{equation*}
\par(2) There exists a constant $C_2$ s.t. for all $\varepsilon \in (0,1)$ and all $t\in [0,T]$,
\begin{equation*}
    \displaystyle\left\|\rho_{g_\varepsilon}\right\|_{\infty} \leq C_{2}.
\end{equation*}   
\end{lemma}
\par To prove Theorem \ref{thm2}, we introduce the following key proposition about the random regularity in the normal regime given in \cite{3}.
\begin{proposition}\label{pro5.5}
  For $d=1$, under Assumptions \ref{ass3.4}-\ref{ass3.5} of 
 $\tilde{a},a_1,a_2,t_0,f^{*}$ and all of its $x,v,z$-derivatives up to total order $K$, the following estimates hold:
\begin{equation}
    \left\|\partial_{z}^{k}\Big[f(x, v, t, z)-f^{*}(x-v t, v, z)\Big]\right\|_{\tilde{a}, t_{0}, 1} \leq C, ~\text{for}~ 0 \leq k \leq K, \label{9}
\end{equation}
and
\begin{equation}
\left\|\partial_{x}^{l} \partial_{z}^{k} E\right\|_{\tilde{a}, t_{0}, l} \leq C,~ \text{for all}~  l \geq 0, ~ k\geq 0, ~ l+k \leq K, \label{10}
\end{equation}
where $f$ is the Landau damping solution to Eq. \eqref{2.1} that satisfies Eq. \eqref{17*} and $E$ is the corresponding electric field,  i.e. $f~\text{and}~E$ form a pair of solutions to Eq. \eqref{2.1}.~Here $C>0$ is a constant.   
\end{proposition}
\par We will prove Theorem \ref{thm2} based on Proposition \ref{pro5.5} in the next section.
\section{Proof of main results}
\label{sec6}
\subsection{Proof of Theorem \ref{thm1}}
\par In the proof of Theorem \ref{thm1}, we first establish an upper bound of
$\sup \limits _{t \in[0, T]} W_{1}(\tilde{f_{\varepsilon}}(t,z), g(t,z))$. After that, each part of this upper bound is estimated separately. 
\par Under the assumption of Theorem \ref{thm1} and by Lemma \ref{lem5.1}, the following inequalities hold for all $z\in A_{z,\delta,M,z_{0}}$:
\begin{equation} 
 W_{1}\left({\tilde f_{\varepsilon}(t,z)}, g(t,z)\right) \leq W_{1}\left(\tilde{f_{\varepsilon}}(t,z), \tilde{g}_{\varepsilon}(t,z)\right)+W_{1}\left(\tilde{g}_{\varepsilon}(t,z), g(t,z)\right), \label{50}
\end{equation}
and
\begin{align}
W_{1}\left(\tilde{f}_{\varepsilon}(t,z), \tilde{g}_{\varepsilon}(t,z)\right)\le\left(1+C_{T}\right) W_{1}\left(f_{\varepsilon}(t,z), g_{\varepsilon}(t,z)\right)\le \left(1+C_{T}\right) W_{2}\left(f_{\varepsilon}(t,z), g_{\varepsilon}(t,z)\right).\label{51}
\end{align}
By Eqs. \eqref{50} and \eqref{51}, the proof of Theorem \ref{thm1} can be established by the following two steps:
\par Step 1:  prove that $\sup \limits _{t \in[0, T]} W_{1}({\tilde g_{\varepsilon}}(t,z), g(t,z))$ uniformly converges to 0 in~$A_{z, \delta,M,z_0}$ as $\varepsilon \to 0$.
\par Step 2:  prove that $ \sup \limits _{t \in[0, T]} W_{2}(f_{\varepsilon}(t,z), g_{\varepsilon}(t,z))$ uniformly converges to 0 in~$A_{z, \delta,M,z_0}$ as $\varepsilon \to 0$.

Now we begin the proof of Theorem \ref{thm1}.

\begin{proof}
\par We divide the proof into two steps.
\par \textbf{Step 1:}
\par Denoting 
\begin{align*}
 G_{\varepsilon}(z)&=\sup \limits_{t\in[0,T]}\left\{\sup \limits _{ ~\theta \in \mathbb{R}^d}\left\|\rho_{\varepsilon}^{\theta}(x,t,z)\right\|_{\infty} 
 \int_{\mathbb{R}^d}\left\|v_{\varepsilon}^{\theta}(x, t,z)+C_{\varepsilon}(x, t,z)-v^{\theta}(x, t,z)\right\|_{\infty} d \mu(\theta) \right.\nonumber\\&\left. + \int_{\mathbb{R}^d}\left\|\rho_{\varepsilon}^{\theta}(x,t,z)-\rho^{\theta}(x,t,z)\right\|_{\infty} d \mu(\theta)\left(\dfrac{1}{2}+\sup \limits _{\theta \in \mathbb{R}^d}\left\|v^{\theta}(x, t,z)\right\|_{\infty}\right)\right\}, d=2,3,
\end{align*}
by Lemma \ref{lem5.2} and the definition of $A_{z, \delta,M,z_0}$, one has such a uniform estimate:
\begin{equation}
     \sup \limits _{t \in[0, T]} W_{1}\left(\tilde{g}_{\varepsilon}, g\right)\le G_{\varepsilon}(z)\le M^2G_{\varepsilon}(z_0), \forall z\in A_{z, \delta,M,z_0},\label{6.3}
\end{equation}
where the first inequality follows by applying Lemma \ref{lem5.2} to each fixed $z$ and the second inequality is achieved from the definition of $A_{z, \delta,M,z_0}$. In fact, by the definition of $A_{z, \delta,M,z_0}$, for any $z\in A_{z, \delta,M,z_0} $, the following inequalities hold:
\begin{align}
  \left\|\rho_{\varepsilon}^{\theta}(x,t,z)\right\|_{\infty}&\le M\left\|\rho_{\varepsilon}^{\theta}(x,t,z_0)\right\|_{\infty}, \nonumber\\ \int_{\mathbb{R}^d}\left\|v_{\varepsilon}^{\theta}(x, t,z)+C_{\varepsilon}(x, t,z)-v^{\theta}(x, t,z)\right\|_{\infty} d \mu(\theta)&\le M\int_{\mathbb{R}^d}\left\|v_{\varepsilon}^{\theta}(x, t,z_0)+C_{\varepsilon}(x, t,z_0)-v^{\theta}(x, t,z_0)\right\|_{\infty} \nonumber \\ &~~~~~~~~~~~~~~~~~~~~~~~~~~~~~~~~~~~~~~~~~~~~~~~~~~~~~~~~~~~~~~~~~~~~d\mu(\theta),\nonumber\\\int_{\mathbb{R}^d}\left\|\rho_{\varepsilon}^{\theta}(x,t,z)-\rho^{\theta}(x,t,z)\right\|_{\infty} d \mu(\theta)&\le M\int_{\mathbb{R}^d}\left\|\rho_{\varepsilon}^{\theta}(x,t,z_0)-\rho^{\theta}(x,t,z_0)\right\|_{\infty} d \mu(\theta),\nonumber\\\left\|v^{\theta}(x, t,z)\right\|&\le M\left\|v^{\theta}(x, t,z_0)\right\|.\label{6.4}
\end{align}
As $M>1$, The last inequality in Eq. \eqref{6.4} indicates 
\begin{align}
    \frac{1}{2}+\left\|v^{\theta}(x, t,z)\right\|\le M\Big(\frac{1}{2}+\left\|v^{\theta}(x, t,z_0)\right\|\Big) \label{6.5}.
\end{align}
By Eq. \eqref{6.5} and the first three inequalities in Eq. \eqref{6.4}, one can derive the second inequality in Eq. \eqref{6.3}.

After taking $z=z_0$, one acquires a deterministic model. By Lemma \ref{lem5.2}, $G_\varepsilon(z_0) \to0,~ \text{as} ~\varepsilon\to 0.$ Thus $\sup \limits _{t \in[0, T]} W_{1}(\tilde{g}_{\varepsilon}, g)$ uniformly converges to 0 in~$A_{z, \delta,M,z_0}$ as $\varepsilon \to 0$. We complete the proof of Step 1.
\par \textbf{Step 2:}
\par  By the assumptions of Theorem \ref{thm1}, Assumption \ref{ass3.2} (i)-(vi)  and Eq. \eqref{22} hold for all $ z\in I_z$. Hence Lemma \ref{lem5.4} holds for all $z\in I_z$.
\par By Lemma \ref{lem5.4},  one yields the estimate
\begin{equation*}
T\le \int_{0}^{T} \tilde{A}(s,z) d s\le T\Biggl(1+\dfrac{1+C_\varepsilon(d)+\sqrt{C_2\max \left \{ C_2,C_\varepsilon (d) \right \} }}{\varepsilon^2}\Biggr),  
\end{equation*}
where for $d=2$,  $ C_{\varepsilon}(d) =C_{\beta_0} / \varepsilon^{2 \max \{\beta_0, \gamma\}}$ for a fixed $\beta_0>2$;  for $d=3$,  $C_{\varepsilon}(d) =C_1 / \varepsilon^{ \max \{38, 3\gamma\}}$.
Thus 
$R_\varepsilon(T,x,z)\le \tilde{R}_\varepsilon(T,x)$,
where 
\begin{align*}  
  \tilde{R} (t,x)=\begin{cases}16 d \exp\Biggl [ {\log \left(\dfrac{x}{16 d}\right) \exp \Biggl(C_{0}t\biggl(1+\dfrac{1+C_\varepsilon(d)+\sqrt{C_2\max \left \{ C_2,C_\varepsilon (d) \right \} }}{\varepsilon^2}\biggr)\Biggr)}\Biggr ], x>16d, \\16d\exp \biggl[\displaystyle\log\Bigl(\dfrac{x}{16d}\Bigr) \exp({C_{0}}t)\biggr]   , 0< x\le 16d.\end{cases}
\end{align*}
\par There exists a function $\phi(\varepsilon)$ such that for all~$0<\varepsilon<\delta$, 
both
$\phi(\varepsilon)\le d $ and $\tilde{R}_\varepsilon(T, \phi(\varepsilon))\le d$  hold, and 
\begin{equation*}
    L(\varepsilon):= \tilde{L}\bigl(\phi(\varepsilon)\bigr) \to 0,~\text{as}~\varepsilon\to 0 ,
\end{equation*}
where
\begin{align*}
\tilde{L} (x)=\begin{cases}16 d \exp\Biggl [ {\log \left(\dfrac{x}{16 d}\right) \exp \Biggl(\dfrac{C_{0} T}{\varepsilon^2}\biggl(1+C_2+C_\varepsilon(d)\biggr)\Biggr)}\Biggr ], x>16d, \\16d\exp \biggl[\displaystyle\log\Bigl(\dfrac{x}{16d}\Bigr) \exp\Bigl(\dfrac{C_{0} T}{\varepsilon^2}\Bigr)\biggr], 0< x\le 16d.\end{cases}    
\end{align*}
For example, let 
\begin{equation*}
 \phi (\varepsilon )=16d\exp \Bigl[-\exp (\varepsilon ^{-K_0} )\Bigr]   
\end{equation*}
with a big enough integer $K_0$.
\par There exists an initial Wasserstein distance $\psi (\varepsilon ,z)\le \phi(\varepsilon),$  for example, one can take
\begin{equation*}
 \psi (\varepsilon ,z)=\dfrac{\phi (\varepsilon )}{1+z^{2} } .  
\end{equation*}
So 
\begin{equation*}
R_\varepsilon\Bigl(T, \psi(\varepsilon,z),z\Bigl)\le \tilde{R}_\varepsilon\Bigl(T, \psi(\varepsilon,z)\Bigl)\le\tilde{R}_\varepsilon\Bigl(T, \phi(\varepsilon)\Bigl)\le d .   
\end{equation*}
As $\psi (\varepsilon ,z) \le \phi(\varepsilon)\le d$, $\psi (\varepsilon ,z)$ satisfies (vi) and (vii) in Assumption 3.1.
\par Then by Lemmas \ref{lem5.3}-\ref{lem5.4}, 
\begin{align}
\sup \limits _{t \in[0, T]} W_{2}\left(f_{\varepsilon}, g_{\varepsilon}\right)&\le16 d \exp \Biggl\{\operatorname { l o g } \biggl( \dfrac { \psi ( \varepsilon,z ) } { 1 6 d } \biggr ) \operatorname {exp} \left[C _ { 0 } T \dfrac { 1 } { \varepsilon ^ { 2 } } \left(1+\left\|\rho_{f_{\varepsilon}}\right\|_{L^{\infty}\left([0, T] ;L_{x}^{\infty}\right)}+\left\|\rho_{g_{\varepsilon}}\right\|_{L^{\infty}\left([0, T] ; L_{x}^{\infty}\right)}\right)\right]\Biggr\}\nonumber\\&\le\tilde{L}\bigl(\psi(\varepsilon,z)\bigl)\le \tilde{L}\bigl(\phi(\varepsilon)\bigl)=L(\varepsilon)\to0.\label{42}   
\end{align}
Eq. \eqref{42}  means $\sup \limits_{t \in[0, T]}W_{2}(f_{\varepsilon}, g_{\varepsilon})$ uniformly converges to 0 in~$A_{z, \delta,M,z_0}$ as $\varepsilon \to 0$, which is the target of Step 2.
\par We complete the proof of Theorem \ref{thm1} through these two steps.
\end{proof}

\subsection{Proof of Theorem \ref{thm2}}
\par In the one-dimensional case, the Vlasov-Poisson system in the quasineutral system Eq. \eqref{7} is reduced to the following form:
\begin{align}
       \left\{\begin{array}{c}\partial_{t} f_\varepsilon(x,v,t,z)+v \partial_{x} f_\varepsilon+E_\varepsilon(x,t,z) \partial_{v} f_\varepsilon=0, ~x\in \mathbb{T},~ t>0, ~v\in \mathbb{R}, \\\varepsilon ^{2} \partial_{x} E_\varepsilon=\displaystyle\int_{-\infty}^{+\infty} f_\varepsilon(x, v, t,z) d v-1,  \\f_\varepsilon(x, v, 0,z)=f_{0, \varepsilon}(x, v,z), \\\displaystyle\int_{0}^{1} \int_{-\infty}^{+\infty} f_\varepsilon(x, v, 0,z) d v d x=1.\end{array}\right.\label{4} 
\end{align}
\par We carry out the proof of Theorem \ref{thm2} in the subsequent two steps.  In Step 1 we prove Eq. \eqref{17} by first performing a variable substitution where $x$ is substituted by $\frac{x}{\varepsilon}$ and $t$ is substituted by $\frac{t}{\varepsilon}$. Subsequently, we estimate
\begin{equation*}
\left\|\partial_{z}^{k}\left[h(x, v, t, z)-f^{*}\Big(\frac{x-v t }{\varepsilon } , v, z\Big)\right]\right\|_{a, t_{0}, 1,m}    
\end{equation*}
based on Eq. \eqref{9} in Proposition \ref{pro5.5} and the relationship between $ \| F  \|  _{a,t_{0},k,m }$ and $ \| F  \|  _{\tilde{a} ,t_{0},k}$. In Step 2 we prove Eq. \eqref{18} and Eq. \eqref{18*} in two cases: when $l\ge1,~l+k\ge2$ and when  $l=0,~k\ge1$. For each case, we first perform the same variable substitution as Step 1, and then the estimation of
\begin{equation*}
\left | \left |  \partial_{x}^{l}\partial_{z}^{k} E_{1}  \right |  \right | _{a,t_{0},l,m }    
\end{equation*}
or
\begin{equation*}
   \left \|  \partial _ {z}^{k} \big[E_{1}(x,t,z)-E_{1}(x_0,t,z)\big]  \right \| _{a,t_{0},0,m}  
\end{equation*}

is derived using Eq. \eqref{10} in Proposition \ref{pro5.5} and the relationship of the two norms.  
Now we begin the proof of Theorem \ref{thm2}.
\begin{proof} 
\par We shall split the proof into two parts.
\par \textbf{Step 1}:
\par Notice that 
\begin{equation}
 \displaystyle\left \| F \right \|  _{a,t_{0},k,m }=\left \| F \right \|  _{\tilde{a} (\varepsilon ),t_{0},k}, \label{37}
\end{equation}  
where~$\tilde{a} (\varepsilon )=a+\frac{1}{\varepsilon ^{m} } $. This equivalence between $ \| F  \|  _{a,t_{0},k,m }$ and $ \| F  \|  _{\tilde{a} (\varepsilon),t_{0},k}$ is significant in applying Proposition \ref{pro5.5} to prove Theorem \ref{thm2} both in this step and the following Step 2. By the assumptions of Theorem \ref{thm2}, for each fixed $\varepsilon\in(0, \delta_1)$, there exist~$\tilde{a} (\varepsilon )=a+\frac{1}{\varepsilon ^{m} },~a_{1}(\varepsilon ),~a_{2}(\varepsilon ),~\varepsilon\text{-independent } t_{0},~\varepsilon$-independent $f^{\ast }$ and all of its $x,v,z$-derivatives up to total order $K$ that satisfy Assumptions \ref{ass3.4}-\ref{ass3.5}. Since for $\varepsilon\in(0, \delta_1)$, by Eq. \eqref{9} in Proposition \ref{pro5.5}, combined with Eq. \eqref{37},  one obtains:
\begin{align*}
\left\|\partial_{z}^{k}\left[f(x, v, t, z)-f^{*}(x-v t, v, z)\right]\right\|_{a, t_{0}, 1,m}=\left\|\partial_{z}^{k}\left[f(x, v, t, z)-f^{*}(x-v t, v, z)\right]\right\|_{\tilde{a}(\varepsilon), t_{0}, 1}\leq C, \nonumber \\ \text{for} ~0 \leq k \leq K, 
\end{align*}
that is,
\begin{equation}
  \sup \limits_{t\ge t_{0}\varepsilon  }\varphi(t, \varepsilon )\left | \left |  \partial_{z}^{k}\left[f(x , v, t, z)-f^{*}(x-v t  , v, z)\right]  \right |  \right | _{L^{\infty } }\le C ,~\text{for} ~0 \leq k \leq K,  \label{49} 
\end{equation}
where
\begin{equation*}
  \displaystyle\varphi  (t, \varepsilon)=\exp\biggl[{\Bigl(a+\frac{1}{\varepsilon ^{m} }\Bigr)t }\biggr]\cdot t^{-1}  .
\end{equation*}
In Eq. \eqref{49},  performing the change of variables $(x, t) \mapsto(\frac{x}{\varepsilon}, \frac{t}{\varepsilon})$, then one has:
\begin{equation}
    \sup \limits_{t\ge t_{0}\varepsilon  }\varphi\Bigl(\frac{t}{\varepsilon }, \varepsilon \Bigr)\left | \left |  \partial_{z}^{k}\left[f\Bigl(\frac{{x}}{\varepsilon }  , v, \frac{{t}}{\varepsilon }, z\Bigr)-f^{*}\Bigl(\frac{x-v t }{\varepsilon } , v, z\Bigr)\right]  \right |  \right | _{L^{\infty } }\le C ,~\text{for} ~0 \leq k \leq K.\label{22*}
\end{equation}
\par Let 
\begin{equation*}
  A(t, \varepsilon)=\displaystyle{\varphi \Big(\frac{t}{\varepsilon } , \varepsilon\Big)}/{\varphi (t, \varepsilon)} =\varepsilon \exp\biggl[{\Big(a+\frac{1}{\varepsilon ^{m} } \Big)\Big(\frac{1}{\varepsilon }-1\Big)t \biggr]},  
\end{equation*}then  
\begin{equation*}
 \displaystyle A(t, \varepsilon)\ge A(t_{0}, \varepsilon )=\varepsilon \exp\biggl[{\Big(a+\frac{1}{\varepsilon ^{m} } \Big)\Big(\frac{1}{\varepsilon }-1\Big)t_0 \biggr]}:=B (\varepsilon )   .
\end{equation*}
By Eq. \eqref{22*}, the following inequality holds:
\begin{align}   
 \sup \limits_{t\ge t_{0}  }&~ \varphi({t, \varepsilon}  )\left | \left |  \partial_{z}^{k}\left[f\Big(\frac{{x}}{\varepsilon }  , v, \frac{{t}}{\varepsilon }, z\Big)-f^{*}\Big(\frac{x-v t }{\varepsilon } , v, z\Big)\right]  \right |  \right | _{L^{\infty } }\nonumber\\ &\le\frac{1}{B(\varepsilon ) } \sup \limits_{t\ge t_{0}  } \varphi\Big(\frac{t}{\varepsilon }, \varepsilon \Big)\left | \left |  \partial_{z}^{k}\left[f\Big(\frac{{x}}{\varepsilon }  , v, \frac{{t}}{\varepsilon }, z\Big)-f^{*}\Big(\frac{x-v t }{\varepsilon } , v, z\Big)\right]  \right |  \right | _{L^{\infty } } \nonumber \\&\le\frac{1}{B(\varepsilon ) } \sup \limits_{t\ge t_{0}\varepsilon  } \varphi\Big(\frac{t}{\varepsilon }, \varepsilon\Big )\left | \left |  \partial_{z}^{k}\left[f\Big(\frac{{x}}{\varepsilon }  , v, \frac{{t}}{\varepsilon }, z\Big)-f^{*}\Big(\frac{x-v t }{\varepsilon } , v, z\Big)\right]  \right |  \right | _{L^{\infty } }\label{23}\\\nonumber &\le \frac{C}{B(\varepsilon ) }=C\exp(at_0)\dfrac{1}{\varepsilon } \exp{\Biggl[- \biggl(\frac{1}{\varepsilon ^{m+1} } -\frac{1}{\varepsilon ^{m}}+\frac{a}{\varepsilon } \biggr)t_{0}\Biggr ] }\to 0,~\text{as}~\varepsilon\to 0, ~ \text{for}~ 0 \leq k \leq K.
\end{align}
 The first ``$\to0$'' in the last line holds since $m\in \mathbb{N}_{+}$, which implies $m+1>1,m+1>m$. This ensures that the dominant term in the exponent is the first term $-\frac{1}{\varepsilon^{m+1}}$.
\par Eq. \eqref{23} indicates that
\begin{equation*}
    \left\|\partial_{z}^{k}\left[h(x, v, t, z)-f^{*}\Big(\frac{x-v t }{\varepsilon } , v, z\Big)\right]\right\|_{a, t_{0}, 1,m} \to 0, ~\text{as}~\varepsilon \to 0,~\text{for}~ 0 \leq k \leq K,
\end{equation*}
which implies the random regularity in the quasineutral regime with a convergence rate of
\begin{equation*}
  \displaystyle O\Biggl(\dfrac{1}{\varepsilon } \exp{\bigg [- \Big(\frac{1}{\varepsilon ^{m+1} } -\frac{1}{\varepsilon ^{m}}+\frac{a}{\varepsilon } \Big)t_{0}\bigg ] }\Biggr)  
\end{equation*}by Eq. \eqref{23}.  We complete the proof of Eq. \eqref{17}.
\par\textbf{ Step 2}:
\par We proceed with Step 2 by considering two cases:
(1) $l\ge 1,~l+k\ge 2$;
(2) $l=0,~k\ge 1$.
\par(1) If~$l\ge 1,~l+k\ge 2$, denote~$l=1+p$, from Eq. \eqref{4} one has  
\begin{equation}
\varepsilon ^{2}  \partial _{x} E_{1}(x,t,z) =\displaystyle\int_{-\infty }^{+\infty } h(x,v,t,z)dv- 1,  \label{41*}  
\end{equation}
namely,
\begin{equation}   
\varepsilon ^{2}  \partial _{x} E_{1}(x,t,z) = \int_{-\infty }^{+\infty } f\Big(\frac{x}{\varepsilon } ,v, \frac{t}{\varepsilon },z\Big)dv- 1.\label{25}
\end{equation}
Then take the $p$-th derivative of $x$ and the $k$-th derivative of $z$ on both sides of Eq. \eqref{25}, which yields:
\begin{equation}
    \varepsilon^{2}\partial _{x}^{l} \partial _ {z}^{k} E_{1}(x,t,z)= \partial _{x}^{p} \partial _{z}^{k} \int_{-\infty }^{+\infty }  f\Big(\frac{x}{\varepsilon } ,v, \frac{t}{\varepsilon },z\Big)dv={\varepsilon} ^{-p}{\partial _{x}^{l} \partial _ {z}^{k} E  }\Big(\frac{x}{\varepsilon } , \frac{t}{\varepsilon },z\Big), \label{55}
\end{equation}
where the second equality comes from taking the $p$-th derivative of $x$ and the $k$-th derivative of $z$ on both sides of $\partial _xE(\frac{x}{\varepsilon}, \frac{t}{\varepsilon},z)=\int_{-\infty }^{+\infty }  f(\frac{x}{\varepsilon},v, \frac{t}{\varepsilon},z) dv-1$.
Eq. \eqref{55} gives
\begin{equation}
  \partial _{x}^{l} \partial _ {z}^{k} E_{1}(x,t,z)= {\varepsilon} ^{-(l+1)}{\partial _{x}^{l} \partial _ {z}^{k} E  }\Big(\frac{x}{\varepsilon } , \frac{t}{\varepsilon },z\Big). \label{55*}
\end{equation}
\par For $\varepsilon\in(0, \delta_1)$,  by Eq. \eqref{10} in Proposition \ref{pro5.5}, combined with Eq. \eqref{37}, one has:
\begin{equation*}
\left\|\partial_{x}^{l} \partial_{z}^{k} E\right\|_{a, t_{0}, l,m}=\left\|\partial_{x}^{l} \partial_{z}^{k} E\right\|_{\tilde{a}(\varepsilon), t_{0}, l} \leq C,    
\end{equation*}
that is,
\begin{equation}
    \sup \limits_{t\ge t_{0}\varepsilon  } r_{p}(t, \varepsilon )\left | \left |  \partial_{x}^{l}\partial_{z}^{k} E(x,t,z)  \right |  \right | _{L^{\infty } }\le C, \label{58}
\end{equation}
where $r_{p}(t, \varepsilon)=\exp[ (a+\frac{1}{\varepsilon})t]\cdot t^{-l}=t^{-p}\varphi(t, \varepsilon)$. Using the change of variables $(x, t) \mapsto(\frac{x}{\varepsilon}, \frac{t}{\varepsilon})$ in Eq. \eqref{58} yields:
\begin{equation}
 \sup \limits_{t\ge t_{0}\varepsilon  } r_{p}\Big(\frac{t}{\varepsilon }, \varepsilon\Big )\left | \left |  \partial_{x}^{l}\partial_{z}^{k} E\Big(\frac{x}{\varepsilon } , \frac{t}{\varepsilon },z\Big)  \right |  \right | _{L^{\infty } }\le C. \label{24}
\end{equation}
Combining Eq. \eqref{55*} and Eq. \eqref{24} gives
\begin{equation*}
   \sup \limits_{t\ge t_{0}\varepsilon  } r_{p}\Big(\frac{t}{\varepsilon }, \varepsilon\Big )\left | \left |  \partial_{x}^{l}\partial_{z}^{k} E_{1}  \right |  \right | _{L^{\infty } }\le C {\varepsilon} ^{-(l+1)} .
\end{equation*}
So
\begin{equation}
     \sup \limits_{t\ge t_{0} } r_{p}\Big(\frac{t}{\varepsilon } , \varepsilon\Big)\left | \left |  \partial_{x}^{l}\partial_{z}^{k} E_{1}  \right |  \right | _{L^{\infty } }\le C {\varepsilon} ^{-(l+1)} .
\end{equation}
\par Notice that
\begin{equation}
r_{p}\Big(\dfrac{t}{\varepsilon }, \varepsilon\Big)/{r_{p}(t, \varepsilon)} =\varepsilon ^{p} A(t, \varepsilon )\ge \varepsilon ^{p}B(\varepsilon ),  \label{6.15}  
\end{equation}
then by Eq. \eqref{6.15} one derives the following estimate:
\begin{align}
    \left | \left |  \partial_{x}^{l}\partial_{z}^{k} E_{1}  \right |  \right | _{a,t_{0},l,m }&=\sup \limits_{t\ge t_{0}  } r_{p}(t, \varepsilon)\left | \left |  \partial_{x}^{l}\partial_{z}^{k} E_{1}  \right |  \right | _{L^{\infty } }\le \frac{1 }{\varepsilon ^{p} B(\varepsilon)} \sup \limits_{t\ge t_{0}  } r_{p}\Big(\frac{t}{\varepsilon }, \varepsilon \Big)\left | \left |  \partial_{x}^{l}\partial_{z}^{k} E_{1}  \right |  \right | _{L^{\infty } }\nonumber\\&\le \dfrac{C {\varepsilon} ^{-(l+1+p)}}{B(\varepsilon )}  =\frac{C {\varepsilon} ^{-2l}}{B(\varepsilon )}\nonumber\\&=C\exp(at_0)\dfrac{1}{\varepsilon ^{2l+1}} \exp{\Biggl[- \biggl(\frac{1}{\varepsilon ^{m+1} } -\frac{1}{\varepsilon ^{m}}+\frac{a}{\varepsilon } \biggr)t_{0}\Biggr ] }.
\end{align}
\par Finally we conclude:
\begin{equation}
   \left | \left |  \partial_{x}^{l}\partial_{z}^{k} E_{1}  \right |  \right | _{a,t_{0},l,m }\to 0,~\text{as}~\varepsilon\to 0,
\end{equation}
with the convergence rate
\begin{equation*}
 \displaystyle O\Biggl(\dfrac{1}{\varepsilon^{2l+1} } \exp{\biggl[- \Bigl(\frac{1}{\varepsilon ^{m+1} } -\frac{1}{\varepsilon ^{m}}+\frac{a}{\varepsilon } \Bigr)t_{0}\biggr ] }\Biggr).   
\end{equation*}
\par (2) If $l=0,~k\ge 1$, for fixed $x_0\in[0,1]$, integrating both sides of  Eq. \eqref{41*} with respect to $x$ from $x_0$ to $x$, it follows
\begin{equation} 
\varepsilon^{2}\big(E_{1}(x, t, z)-E_{1}\left(x_{0}, t, z\right)\big)=\int_{x_{0}}^{x} \int_{-\infty}^{+\infty} h(s, v, t, z) d v d s-x+x_{0}.\label{6.38}
\end{equation}
Integrating both sides of the second line of Eq. \eqref{2.1} in the one-dimensional case with respect to $x$ from $\frac{x_0}{\varepsilon}$ to $x$, one has
 \begin{equation}
   E(x, t, z)-E\left(\frac{x_{0}}{\varepsilon}, t, z\right)=\int_{\frac{x_{0}}{\varepsilon}}^{x} \int_{-\infty}^{+\infty} f(s, v, t, z) d v d s-x+\frac{x_{0}}{\varepsilon}.\label{6.18}
\end{equation}
So taking the $k$-th derivative of $z$ on both sides of Eq. \eqref{6.38} yields the following equality:
\begin{align}
  \varepsilon^{2}\partial^{k}  _{z}\big[E_{1}(x, t, z)-E_{1}\left(x_{0}, t, z\right)\big]&=\displaystyle\partial^{k}  _{z}\int_{x_{0}}^{x}  \int_{-\infty }^{+\infty } h(s,v,t,z)dvds=\partial^{k}  _{z} \int_{x_{0}}^{x}  \int_{-\infty }^{+\infty } f\Big(\dfrac{s}{\varepsilon }   ,v, \frac{t}{\varepsilon },z\Big)dvds\nonumber\\ &=\varepsilon \partial^{k}  _{z}\int_{\frac{x_{0}}{\varepsilon }}^{\frac{x}{\varepsilon } }  \int_{-\infty }^{+\infty } f\Big(s,v , \frac{t}{\varepsilon},z\Big)dvds= \varepsilon   \partial^{k}  _{z}\bigg[ E(\frac{x}{\varepsilon}, \frac{t}{\varepsilon}, z)-E\left(\frac{x_{0}}{\varepsilon}, \frac{t}{\varepsilon}, z\right)\bigg],\label{6.19*}
\end{align}
where the second equality follows from Eq. \eqref{3.7}, and the last equality is derived by taking the $k$-th derivative of $z$ on both sides of Eq. \eqref{6.18}.
\par By Eq. \eqref{10} in proposition \ref{pro5.5} and Eq. \eqref{37}, we obtain
\begin{align}
 \left\| \partial_{z}^{k}\bigg[E(x,t,z)-E\Big(\frac{x_0}{\varepsilon},t,z\Big)\bigg]\right\|_{a, t_{0}, 0,m}&=\left\|\partial_{z}^{k}\bigg[E(x,t,z)-E\Big(\frac{x_0}{\varepsilon},t,z\Big)\bigg]\right\|_{\tilde{a}(\varepsilon), t_{0}, 0}\nonumber\\& \leq  \left\| \partial_{z}^{k}E(x,t,z)\right\| _{\tilde{a}(\varepsilon), t_{0}, 0}+\left\| \partial_{z}^{k}E\Big(\frac{x_0}{\varepsilon},t,z\Big)\right\| _{\tilde{a}(\varepsilon), t_{0}, 0}\nonumber\\&\le2\left\| \partial_{z}^{k}E(x,t,z)\right\| _{\tilde{a}(\varepsilon), t_{0}, 0}\le2C.\label{6.19}
\end{align}
Then, defining $\tilde{E}(x,t,z)=E(x,t,z)-E\Big(\frac{x_0}{\varepsilon},t,z\Big)$, Eq. \eqref{6.19} gives
\begin{equation}
   \sup \limits_{t\ge t_{0}  } r_{-1}(t, \varepsilon )\left | \left |  \partial_{z}^{k} \tilde E(x,t,z)  \right |  \right | _{L^{\infty } }\le 2C.\label{6.20}
\end{equation}
Here $p=-1$ since $l=0$. Performing the change of variables $(x, t) \mapsto(\frac{x}{\varepsilon}, \frac{t}{\varepsilon})$ in Eq. \eqref{6.20} yields:
\begin{equation}
 \sup \limits_{t\ge t_{0}\varepsilon  } r_{-1}\Big(\frac{t}{\varepsilon }, \varepsilon\Big )\left | \left |  \partial_{z}^{k} \tilde E\Big(\frac{x}{\varepsilon } , \frac{t}{\varepsilon },z\Big)  \right |  \right | _{L^{\infty } }\le 2C. \label{6.21}
 \end{equation}
Combining Eq. \eqref{6.19*} and Eq. \eqref{6.21}, one obtains
\begin{align*}
  \sup \limits_{t\ge t_{0} } r_{-1}\Big(\frac{t}{\varepsilon } , \varepsilon\Big)\left | \left | \partial_{z}^{k}\left[E_{1}(x, t, z)-E_{1}\left(x_{0}, t, z\right)\right]  \right |  \right | _{L^{\infty } }&\le \sup \limits_{t\ge t_{0}\varepsilon  } r_{-1}\Big(\frac{t}{\varepsilon }, \varepsilon\Big )\left | \left |  \partial_{z}^{k}\left[E_{1}(x, t, z)-E_{1}\left(x_{0}, t, z\right)\right] \right |  \right | _{L^{\infty } }\\&\le{\varepsilon} ^{-1} \sup \limits_{t\ge t_{0}\varepsilon  } r_{-1}\Big(\frac{t}{\varepsilon }, \varepsilon\Big )\left | \left |  \partial^{k}  _{z}\bigg[ E(\frac{x}{\varepsilon}, \frac{t}{\varepsilon}, z)-E\left(\frac{x_{0}}{\varepsilon}, \frac{t}{\varepsilon}, z\right)\bigg] \right |  \right | _{L^{\infty } } \\&={\varepsilon} ^{-1} \sup \limits_{t\ge t_{0}\varepsilon  } r_{-1}\Big(\frac{t}{\varepsilon }, \varepsilon\Big )\left | \left |  \partial_{z}^{k} \tilde E\Big(\frac{x}{\varepsilon } , \frac{t}{\varepsilon },z\Big)  \right |  \right | _{L^{\infty } }\\&\le 2C {\varepsilon} ^{-1} .
\end{align*}
Then by Eq. \eqref{6.15} one acquires the following estimate:
\begin{align*}
    \left | \left |  \partial_{z}^{k}\left[E_{1}(x, t, z)-E_{1}\left(x_{0}, t, z\right)\right]  \right |  \right | _{a,t_{0},0,m }&=\sup \limits_{t\ge t_{0}  } r_{-1}(t, \varepsilon)\left | \left |  \partial_{z}^{k}\left[E_{1}(x, t, z)-E_{1}\left(x_{0}, t, z\right)\right]  \right |  \right | _{L^{\infty } }\\&\le \frac{\varepsilon}{B(\varepsilon)} \sup \limits_{t\ge t_{0}  } r_{-1}\Big(\frac{t}{\varepsilon }, \varepsilon \Big)\left | \left |  \partial_{z}^{k}\left[E_{1}(x, t, z)-E_{1}\left(x_{0}, t, z\right)\right] \right |  \right | _{L^{\infty } }\nonumber\\&\le\frac{2C}{B(\varepsilon )}\nonumber=2C\exp(at_0)\dfrac{1}{\varepsilon } \exp{\Biggl[- \biggl(\frac{1}{\varepsilon ^{m+1} } -\frac{1}{\varepsilon ^{m}}+\frac{a}{\varepsilon } \biggr)t_{0}\Biggr ] }.
\end{align*}
So
\begin{align*}
 {\left \| \partial_{z}^{k}\left[E_{1}(x, t, z)-E_{1}\left(x_{0}, t, z\right)\right] \right \| }_{a,t_{0}, 0,m} \to 0,~\text{as}~\varepsilon\to 0. 
\end{align*}

By the arbitrariness of $x_0$, Eq. \eqref{18*} holds for any $x_0\in[0,1]$.
Moreover, the convergence rate exhibits 
\begin{equation*}
 \displaystyle O\Biggl(\dfrac{1}{\varepsilon } \exp{\bigg[- \Bigl(\frac{1}{\varepsilon ^{m+1} } -\frac{1}{\varepsilon ^{m}}+\frac{a}{\varepsilon } \Bigr)t_{0}\bigg ] }\Biggr).   
\end{equation*}
\par The proofs of Eq. \eqref{18} and Eq. \eqref{18*} are accomplished by discussing these two cases.
\par We complete the proof of Theorem \ref{thm2} by proving Eq. \eqref{17}, Eq. \eqref{18}, and Eq. \eqref{18*} in turn.
\end{proof}
\section{Conclusion}
\label{sec7}
\par In this paper, we prove the uniform convergence with respect to $z$ 
in the Wasserstein distance between the Vlasov-Poisson system with initial uncertainty in the quasineutral regime and its quasineutral limit system for the two-dimensional and three-dimensional cases. This is achieved by imposing appropriate constraints on $z$, establishing an upper bound, and estimating each term of the bound. In this way, we overcome the analytical challenges arising from the randomness. The result shows that the convergence in the Wasserstein distance is insensitive to the uncertainty in the initial data.
 Furthermore, we prove the random regularity of the one-dimensional Vlasov-Poisson system with initial uncertainty in the quasineutral regime, with the random regularity for higher dimensions reserved for future work, which would pose additional difficulties in the analysis. By constructing a new norm that incorporates the quasineutral parameter $\varepsilon$, we derive the $z$-regularity of the distribution function and the electric field of the uncertain Vlasov-Poisson system in the quasineutral regime by this norm, with the convergence rate explicit in the quasineutral parameter. These results firstly quantify the propagation of the initial uncertainty
 for the uncertain Vlasov-Poisson system in the quasineutral regime, which is critical for developing numerical schemes for multiscale kinetic equations with uncertainties, laying theoretical foundations for critical technologies in engineering such as controlling the impact of initial disturbances in semiconductor manufacturing 
 techonologies.
\par This paper primarily focuses on the theoretical properties of the Vlasov-Poisson system with random initial inputs in the quasineutral regime. These theoretical results help us understand how the uncertainty propagates, and how the uncertainty affects the solution for large time. In the future, one may verify the results of this study through numerical experiments, and develop efficient numerical algorithms based on these findings, leveraging numerical tools such as polynomial chaos expansion. Additionally, more factors such as magnetic field and the high-field limit can be taken into consideration. Moreover, the uncertainty quantification of  the Vlasov-Poisson system with more general Landau damping and other kinetic equations with similar damping phenomena are also valuable and promising directions for further investigation.
\section*{Acknowledgments}
This work was supported by National Key R\&D Program of China (2024YFA1016100), National Natural Science Foundation of China (12531016), and College Students' Innovative Entrepreneurial Training Plan Program (IPP30115). W. Wang wishes to express sincere gratitude to Y. Lin for providing careful guidance throughout the research process. The authors would like to thank Y. Zhu and H. Heng for their insightful discussions during the early stages of this research.


\end{document}